\documentclass{amsart}

\usepackage{amssymb} \usepackage{amsfonts} \usepackage{amsmath}
\usepackage{amsthm} \usepackage{epsfig} 
\usepackage{color}
\usepackage{amscd}
\usepackage[all]{xy}

\theoremstyle{definition}

\theoremstyle{remark}

\newcommand{\matr} [4] {\big({\tiny\begin{array}{@{}c@{\ }c@{}} #1 & #2 \\ #3 & #4 \\ \end{array}} \big)}

\newcommand{\Iso}{{\rm Isom}}

\newcommand{\matR} {\ensuremath {\mathbb{R}}}

\newcommand{\matZ} {\ensuremath {\mathbb{Z}}}

\newcommand{\matH} {\ensuremath {\mathbb{H}}}

\newcommand{\matRP} {\ensuremath {\mathbb{RP}}}

\newcommand{\Vol} {\ensuremath {{\rm Vol}}}

\author{Bruno Martelli}
\address{Dipartimento di Matematica ``Tonelli'', Largo Pontecorvo 5, 56127 Pisa, Italy}
\email{martelli at dm dot unipi dot it}

\title{Hyperbolic four-manifolds}

\begin{document}

\begin{abstract}
This is a short survey on finite-volume hyperbolic four-manifolds. We first describe some general theorems, and then focus on the geometry of the concrete examples that we found in the literature. 

The starting point of most constructions is an explicit reflection group $\Gamma$ acting on $\matH^4$, together with its Coxeter polytope $P$. Hyperbolic manifolds then arise either algebraically from the determination of torsion-free subgroups of $\Gamma$, or more geometrically by assembling copies of $P$. 

We end the survey by raising a few open questions.
\end{abstract}

\maketitle

\section*{Introduction}

Hyperbolic $n$-manifolds exist for every $n\geq 2$, but the small dimensions $n=2$ and $n=3$ have been the object of a much wider and more intense study than others. The reason for that is of course well-known: thanks to Riemann's uniformization theorem and Thurston's geometrization, ``most'' closed manifolds in dimension 2 and 3 are hyperbolic, so hyperbolic geometry is the most powerful tool available to understand the topology of manifolds of dimension 2 and 3.

The role of hyperbolic geometry in dimension $n=4$ is less clear; as far as we know (which is not much), hyperbolic four-manifolds seem quite sporadic in the enormous and wild world of smooth four-manifolds. The reason is at least threefold: there are really many smooth four-manifolds around, there is no canonical decomposition available whatsoever (in contrast to dimension three), and there are only ``few'' hyperbolic four-manifolds because hyperbolic geometry is more rigid in dimension $n \geq 4$ than it is in dimensions $n=2$ and $3$. 

The fact that there are ``fewer'' hyperbolic four-manifolds than three-manifolds is of course debatable. The main gift of the three-dimensional hyperbolic world, which lacks in higher dimension, is of course the hyperbolic Dehn filling theorem: a notable consequence is that there are infinitely many closed hyperbolic three-manifolds with volume smaller than 3, while there are only finitely many complete hyperbolic four-manifolds with volume smaller than any number.

However, this finite number of hyperbolic manifolds with bounded volume can be very big and is still completely unknown (we only know that it grows roughly factorially with the bound \cite{BGLM}). In dimension four, the Gauss-Bonnet Theorem says that 
\begin{equation} \label{Vol:chi:eqn}
\Vol(M) = \frac 43 \pi^2 \chi(M)
\end{equation}
on every finite-volume complete hyperbolic four-manifold $M$, hence volume and Euler characteristic are roughly the same thing. 

One of the most important papers on hyperbolic four-manifolds is the Ratcliffe-Tschantz census \cite{RT0} that tabulates 1171 cusped manifolds having $\chi=1$, that is with smallest possible volume. All these manifolds are constructed by some particularly simple side-pairings from a single polytope, the ideal 24-cell. We know absolutely nothing about the actual number of hyperbolic four-manifolds with $\chi=1$, which could be \emph{a priori} much bigger than 1171.

In this paper we survey all the concrete examples of finite-volume hyperbolic four-manifolds that we were able to find in the literature. The starting point of most constructions is an explicit reflection group $\Gamma$ acting on $\matH^4$ together with its corresponding Coxeter polytope $P$. Hyperbolic manifolds then arise either algebraically from the determination of torsion-free subgroups of $\Gamma$, or more geometrically by assembling copies of $P$. Although simple in theory, both strategies  suffer from the combinatorial complexities of both $\Gamma$ and $P$, indeed the determination of torsion-free subgroups has been obtained only in some cases (using computers), and only few very symmetric polytopes $P$ could be assembled successfully.

The paper is organized as follows. We start by introducing reflection groups and Coxeter polytopes in Section \ref{polytopes:section}. We then turn to hyperbolic manifolds in Sections \ref{manifolds:section} and \ref{manipulating:section}. (The former deals with the crude constructions, the latter with some more refined aspects like Dehn fillings, isometries, etc.) Finally, we raise a few open questions in Section \ref{questions:section}.

\subsection*{Acknowledgements}
The author warmly thank Mikhail Belolipetsky, Vincent Emery, Ruth Kellerhals, Alexander Kolpakov, John Ratcliffe, Stefano Riolo, and Steven Tschantz for many fruitful discussions. He also thanks the AIM and its SQuaRE program for creating a very nice and stimulating environment for discussion and research.

\section{Four-dimensional hyperbolic polytopes} \label{polytopes:section}
We list here some four-dimensional hyperbolic polytopes that have been considered by various authors, mostly with the aim of constructing hyperbolic four-manifolds. 

In dimension three, using Thurson's equations one can construct plenty of cusped hyperbolic manifolds by assembling simplices of very different kinds. In contrast, in dimension four the only polytopes that people have been able to use to construct manifolds are essentially the \emph{Coxeter polytopes}, namely those with dihedral angles $\frac \pi k$, and some few more with angles $\frac {2\pi} k$. Coxeter polytopes are particularly nice because they are the fundamental domains of discrete \emph{reflection groups}, that is discrete groups $\Gamma$ of isometries in $\matH^4$ generated by reflections.

\subsection{Coxeter polytopes}
We briefly introduce the theory of hyperbolic Coxeter polytopes. For a more general reference to the subject, see for instance \cite{H, R}.

A \emph{finite polytope} $P$ in $\overline{\matH^n}$ is the convex hull of finitely many points, and we suppose for simplicity that $P$ has non-empty interior. The boundary of $P$ stratifies into faces of various dimensions, named \emph{vertices}, \emph{edges}, $\ldots$, and \emph{facets}. Each vertex is \emph{finite} or \emph{ideal} according to whether it lies in $\matH^n$ or $\partial \matH^n$.

A \emph{Coxeter polytope} is a finite polytope $P$ whose dihedral angles divide $\pi$. The adjacency of facets in $P$ is usually encoded via a graph, in which every node represents a facet of $P$ and decorated edges describe the way two distinct facets meet (or do not meet). 

When $P$ is combinatorially a simplex, every pair of facets meet at some angle $\frac \pi k$, and one labels the corresponding edge with the number $k$. Since one commonly gets many right angles it is customary to draw only the edges with $k\geq 3$, labeled with $k$. Moreover, an unlabeled edge is tacitly assumed to have $k=3$.

For more complicated polytopes, we use the same rules, plus a label $\infty$ on edges to indicate two facets that do not meet. One can also use some more refined label that encodes whether the two hyperplanes containing the facets are asymptotically parallel or ultraparallel, and their mutual distance. Of course the graph can become quite complicated if the polytope has many facets (like for instance in a dodecahedron).

Reflections along the hyperplanes containing the facets of $P$ generate a discrete group $\Gamma < \Iso(\matH^n)$ of isometries with fundamental domain precisely $P$. This means that the polytopes $\{\gamma P\ |\ \gamma \in \Gamma\}$ form a tessellation of $\matH^n$, i.e.~they cover $\matH^n$ and intersect only in common faces. A presentation for $\Gamma$ is
$$\langle R_i \ |\ (R_iR_j)^{k_{ij}}\rangle$$
where $R_i$ is the reflection along the hyperplane containing the facet $f_i$, we set $k_{ii} =1$, and two distinct facets $f_i, f_j$ intersect at an angle $\frac{\pi}{k_{ij}}$. (When they do not intersect we assume $k_{ij}=\infty$ and the relation is omitted.)

Every Coxeter polytope $P$ may be interpreted as an orbifold and has a rational Euler characteristic $\chi(P)$ which may be computed as
$$\chi(P) = \sum_f \frac{(-1)^{\dim f}}{|{\rm Stab} f|}$$
where the sum is over all faces $f$ of $P$ and Stab$f< \Gamma$ is the stabilizer of $f$, see for instance \cite{Br}.
The formula $\Vol(P) = \frac 43 \pi^2 \chi(P)$ holds also here in dimension $4$.

By Selberg's Lemma, there are plenty of torsion-free finite-index subgroups $\Gamma'<\Gamma$ defining hyperbolic manifolds $M=\matH^n/_{\Gamma'}$ that orbifold-cover $P$ and is tessellated into finitely many copies of $P$. A natural two-steps method to construct hyperbolic manifolds is to provide a Coxeter polytope $P$, and then a finite-index torsion-free subgroup $\Gamma' < \Gamma$.

\subsection{Hyperbolic simplices}
\begin{figure}
 \begin{center}
  \includegraphics[width = 9 cm]{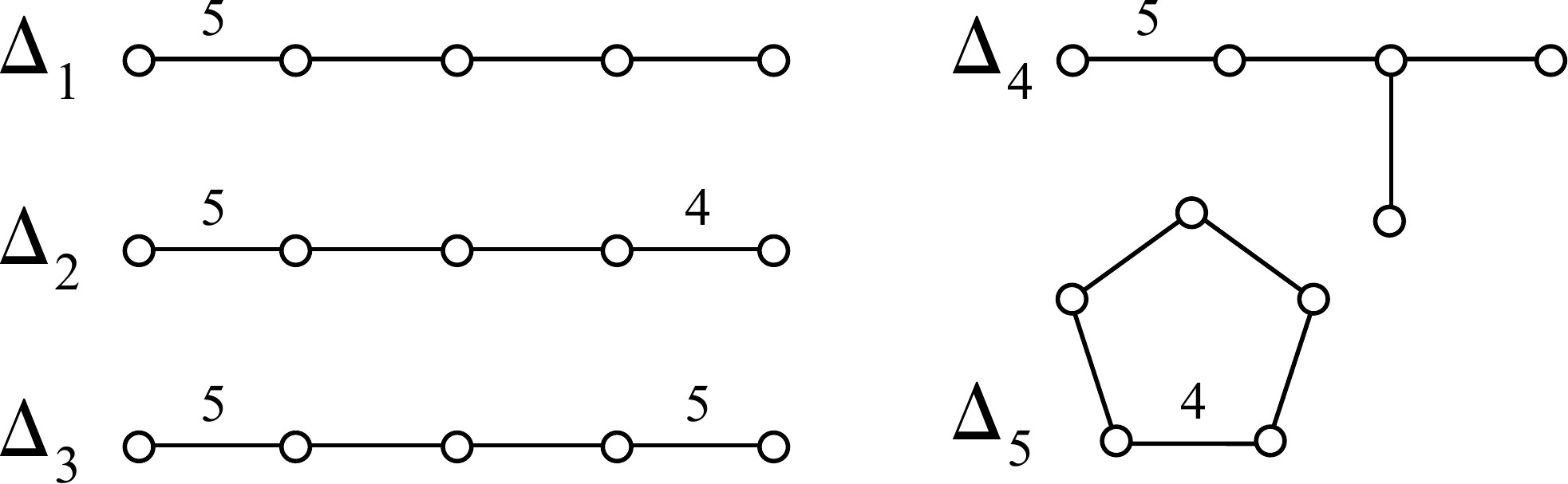}
 \end{center}
 \caption{The five compact Coxeter simplices in $\matH^4$.}  \label{simplexes:fig}
\end{figure}

The first Coxeter polytopes one investigates are of course the simplices, and these were classified by Vinberg in \cite{Vin}. There are five compact Coxeter simplices $\Delta_1,\ldots, \Delta_5$ in $\matH^4$ and they are shown in Figure \ref{simplexes:fig}. Each $\Delta_i$ is the fundamental domain of a cocompact discrete group $\Gamma_i$ generated by reflections along its facets. Their Euler characteristics are quite small:
$$\chi(\Delta_1) = \frac 1{14400},\ \chi(\Delta_2) = \frac{17}{28800},\ \chi(\Delta_3) = \frac{13}{7200},\ \chi(\Delta_4) = \frac{17}{14400},\ \chi(\Delta_5) = \frac{11}{5760}.$$

The five groups $\Gamma_1,\ldots, \Gamma_5$ are all arithmetic, and the first four $\Gamma_1,\ldots, \Gamma_4$ are commensurable. In particular $\Gamma_4 < \Gamma_2$ with index two, and this is the unique direct inclusion between these groups \cite[Theorem 4]{JKRT2}. 
The simplex $\Delta_1$ is the smallest arithmetic four-dimensional hyperbolic orbifold \cite{Be}.

\begin{figure}
 \begin{center}
  \includegraphics[width = 9 cm]{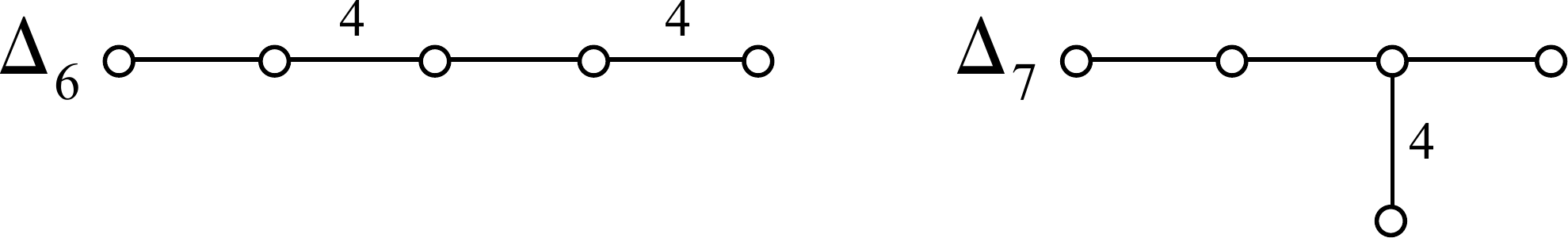}
 \end{center}
 \caption{Two (among the nine) non-compact Coxeter simplices in $\matH^4$.}  \label{ideal_simplexes:fig}
\end{figure}

There are also nine non-compact Coxeter simplices (which have some ideal vertices) and they are all commensurable and arithmetic. Two of them $\Delta_6$ and $\Delta_7$ are shown in Figure \ref{ideal_simplexes:fig} and have
$$\chi(\Delta_6) = \frac 1{1152}, \quad \chi(\Delta_7) = \frac 1{1920}.$$
Let $\Gamma_6$ and $\Gamma_7$ be the corresponding Coxeter groups. Every other Coxeter group arising from the remaining seven simplices is a finite-index subgroup of one of these two (sometimes both), see \cite[Theorem 4]{JKRT2}. 

The simplex $\Delta_7$ is the smallest non-compact arithmetic four-dimensional hyperbolic orbifold \cite{HK}.

\subsection{Regular polytopes}
Luckily, there is no shortage of regular polytopes in dimension four: there are as much as six of them and they are listed in Table \ref{regular:table}. 

By shrinking or inflating a regular polytope $P$ in $\matH^4$ we get a family of hyperbolic regular polytopes with dihedral angles that vary on a segment $[\theta_0,\theta_1) \subset [0,2 \pi]$, where $\theta_1$ is the dihedral angle of $P$ in its euclidean version and $\theta_0$ is the dihedral angle of $P$ in its ideal hyperbolic version. The angle $\theta_0$ is in turn the dihedral angle of the (regular) euclidean vertex figure, which lies in a (euclidean) horosphere. Both angles $\theta_0$ and $\theta_1$ can be spotted from the beautiful book of Coxeter \cite{Cox}, see also \cite{Cox2}.

\begin{table}
\begin{center}
\begin{tabular}{c||ccccc||c}
\phantom{\Big|} \!\! name & vertices & edges & faces & facets & vertex figure & Schl\"afli\\
 \hline \hline
\phantom{\Big|} \!\!  simplex & 5 & 10 & 10 & 5 tetrahedra & tetrahedron & $\{3,3,3\}$\\
\hline
\phantom{\Big|} \!\!  hypercube & 16 & 32 & 24 & 8 cubes & tetrahedron & $\{4,3,3\}$ \\
\phantom{\Big|} \!\!  $16$-cell & 8 & 24 & 32 & 16 tetrahedra & octahedron & $\{3,3,4\}$ \\
\hline
\phantom{\Big|} \!\!  $24$-cell & 24 & 96 & 96 & 24 octahedra & cube & $\{3,4,3\}$ \\
\hline
\phantom{\Big|} \!\!  $120$-cell & 600 & 1200 & 720 & 120 dodecahedra & tetrahedron & $\{5,3,3\}$ \\
\phantom{\Big|} \!\!  $600$-cell & 120 & 720 & 1200 & 600 tetrahedra & icosahedron & $\{3,3,5\}$
\end{tabular}
\vspace{.2 cm}
\caption{The six regular polytopes in dimension four. The dual pairs are grouped (the simplex and the 24-cell are self-dual) and the Schl\"afli symbol is shown.
}
\label{regular:table}
\end{center}
\end{table}

For instance, for the 24-cell we have $[\theta_0, \theta_1) = [\frac{\pi}2, \frac {2\pi}3).$
This wonderfully symmetric polytope has dihedral angle $\frac {2\pi}3$ in the euclidean version, and $\frac{\pi}2$ in the ideal hyperbolic version: the two versions tessellate $\matR^4$ and $\matH^4$ respectively. The 24-cell is the unique regular polytope (in all dimensions $n \geq 3$) whose euclidean and ideal hyperbolic versions are both tessellating.

The generous 120-cell furnishes a big segment of angles, expressed in degrees as $[\theta_0, \theta_1) = [70^\circ 32', 144^\circ)$, which contains $72^\circ$, $90^\circ$, $120^\circ$, that is $\frac{2\pi}5, \frac \pi 2, \frac{2\pi}3$. The other regular hyperbolic polytopes with dihedral angles dividing $2\pi$ are the simplex and the hypercube, both with angles $\frac{2\pi}5$. All these polytopes tessellate $\matH^4$. (Note that those with angles $\frac{2\pi}3$ and $\frac{2\pi}5$ are not Coxeter polytopes.)

Summing up, there are six regular hyperbolic polytopes tessellating $\matH^4$, one ideal and five compact. The tessellations are listed in Table \ref{tessellations:table}.

Let $P$ be a regular $n$-dimensional polytope. A \emph{flag} is a selection of an $i$-dimensional face $f_i$ for each $i=0,\ldots, n$ such that $f_i\subset f_{i+1}$. The barycenters $v_i$ of $f_i$ span a simplex $\Delta\subset P$ called the \emph{characteristic simplex} of $P$. Since the isometry group $\Gamma$ of $P$ acts freely and transitively on flags, the characteristic simplex $\Delta$ is a fundamental domain for $\Gamma$. The characteristic simplices of the tessellating regular polytopes in Table \ref{tessellations:table} are the Coxeter simplices $\Delta_1, \Delta_2, \Delta_3, \Delta_6$ encountered in Figure \ref{simplexes:fig} and \ref{ideal_simplexes:fig}.

\begin{table}
\begin{center}
\begin{tabular}{c||cccc||c}
\phantom{\Big|} \!\! polytope & dihedral angle & vertex figure & char simplex & $\chi$ & Schl\"afli\\
 \hline \hline
\phantom{\Big|} \!\!  simplex & $\frac {2\pi}5$ & $600$-cell & $\Delta_1$ & $\frac 1 {120}$ & $\{3,3,3,5\}$\\
\phantom{\Big|} \!\!  120-cell & $\frac {2\pi}3$ & simplex & $\Delta_1$ & 1 & $\{5,3,3,3\}$\\
 \hline
\phantom{\Big|} \!\!  hypercube & $\frac {2\pi}5$ & 600-cell & $\Delta_2$ & $\frac{17}{75}$ & $\{4,3,3,5\}$ \\
\phantom{\Big|} \!\!  120-cell & $\frac \pi 2$ & 16-cell & $\Delta_2$ & $\frac{17}2$ & $\{5,3,3,4\}$ \\
\hline
\phantom{\Big|} \!\!  120-cell & $\frac {2\pi} 5$ & 600-cell & $\Delta_3$ & $26$ & $\{5,3,3,5\}$ \\
\hline
\phantom{\Big|} \!\!  24-cell & $\frac \pi 2$ & & $\Delta_6$ & 1 & $\{4,3,4,3\}$ 
\end{tabular}
\vspace{.2 cm}
\caption{The six hyperbolic tessellations of $\matH^4$ via regular polytopes. The dual pairs are grouped, and the characteristic simplex (from Figure \ref{simplexes:fig}), Euler characteristic $\chi$, and Schl\"afli symbol are shown. The 24-cell is ideal and hence there is no vertex figure there. (In a sense, the vertex figure is the tessellation $\{3,4,3\}$ of $\matR^3$ into cubes. Following this line, it also makes sense to define a dual tessellation with Schl\"afli symbol $\{3,4,3,4\}$ with infinite polyhedra.)
}
\label{tessellations:table}
\end{center}
\end{table}

Every finite-volume hyperbolic four-manifold $M$ that is tessellated into some hyperbolic regular polytopes is an orbifold cover of the corresponding characteristic simplex $\Delta_1, \Delta_2, \Delta_3$, or $\Delta_6$. All compact manifolds of this kind are commensurable, since $\Delta_1, \Delta_2$, and $\Delta_3$ are. Hence we get overall two commensurability classes, one compact and one non-compact. A closed hyperbolic four-manifold belonging to one of these commensurability classes may not tesselate into regular polytopes.

\subsection{More right-angled polytopes} \label{more:right:subsection}
There is today no complete classification of hyperbolic Coxeter polytopes in dimension four. Note that every right-angled polytope $P$ can be mirrored along any of its facets to give a new right-angled polytope with twice its volume: thus starting from the right-angled 24- or 120-cell we can construct infinitely many complicated right-angled Coxeter polytopes, with arbitrary big volume and number of facets.

There are also right-angled polytopes containing both ideal and finite vertices. One such polyhedron $P_4$ was constructed in \cite{RT0} by cutting the ideal 24-cell along the 4 coordinate hyperplanes into 16 isometric pieces. More precisely, the construction goes as follows: consider the regular ideal 24-cell in the Poincar\'e disc model $D^4$, as the convex hull of the points
$$(\pm 1, 0,0,0), \ (0,\pm 1,0,0),\ ( 0, 0,\pm 1,0),\ (0, 0,0,\pm 1),\ \big(\!\pm\! \tfrac 12, \pm \tfrac 12, \pm \tfrac 12, \pm \tfrac 12\big).$$
The polytope $P_4$ is the intersection of the ideal 24-cell with the positive hexadecant of $\matR^4$. It has five finite vertices
$$(0,0,0,0),\ \big(0,\tfrac 13, \tfrac 13, \tfrac 13\big), \ \big(\tfrac 13, 0, \tfrac 13, \tfrac 13\big) \ \big(\tfrac 13, \tfrac 13, 0, \tfrac 13\big), \big(\tfrac 13, \tfrac 13, \tfrac 13, 0 \big) $$
and five ideal vertices
$$(1,0,0,0), \ (0,1,0,0), \ (0,0,1,0), \ (0,0,0,1), \ \big(\tfrac 12, \tfrac 12, \tfrac 12, \tfrac 12\big).$$
The latter are the ideal vertices of a regular ideal simplex $S$. Indeed $P_4$ is obtained by attaching to $S$ five simplices, each with one finite vertex and four ideal ones. Each of the four exterior facets of the attached simplices matches with an adjacent one, hence $P_4$ has 10 facets overall. Each facet of $P_4$ is a hyperbolic right-angled polyhedron $P_3$ with three ideal vertices and two finite ones. 

The polyhedra $P_3$ and $P_4$ are the first members of a family $P_3, \ldots, P_8$ of right-angled hyperbolic polytopes of dimension $3,\ldots, 8$ that have both finite and ideal vertices, see \cite{PV}. There are no finite-volume right-angled polytopes in dimension $n\geq 13$, see \cite{PV, D}.

Of course the orbifold $P_4$ is commensurable with the ideal 24-cell. It turns out that, among ideal right-angled polytopes, the 24-cell is the unique one having smallest volume, and also the smallest number of facets \cite{K}.

\subsection{Uniform polytopes} \label{rectified:subsection}
A polytope is \emph{uniform} if its isometry group acts transitively on the vertices. For instance, the convex hull of the centers of the $k$-faces of a regular polytope (for some fixed $k$) is a uniform polytope.

A natural example is the \emph{rectified simplex}, which is the convex hull of the midpoints of the edges of a simplex $S$. In dimension three, a rectified simplex is a regular octahedron. 
In higher dimensions, the rectified simplex is not regular. 

In dimension four, the rectified simplex has 10 vertices and we may describe it in $\matR^5$ as the convex hull of $(1,1,1,0,0)$ and all the other points obtained by permuting these coordinates. It has 10 facets: five regular tetrahedra (obtained by truncating the vertices of $S$) and five regular octahedra (obtained by rectifying the facets of $S$).

The hyperbolic version which is relevant for us is the \emph{ideal hyperbolic rectified simplex} $P$, with all 10 vertices in $\partial \matH^4$. The symmetries of $P$ force the vertex figure to be a euclidean prism made of two horizontal equilateral triangles and three vertical squares. The dihedral angles of this prism are $\frac \pi 2$ and $\frac \pi 3$, and hence these are also the dihedral angles of $P$ which is therefore a Coxeter polytope. The facets of $P$ are ideal regular tetrahedra and octahedra: every tetrahedron meets four octahedra with angles $\frac \pi 2$, and two adjacent octahedra meet at angles $\frac \pi 3$. We have $\chi(P) = \frac 16$, see \cite{KS}.

This beautiful hyperbolic polytope has been used to construct many manifolds in \cite{KS} and then also in \cite{S2}.

\subsection{Deforming Coxeter polytopes}
Finite-volume orbifold groups are rigid in dimension four, but infinite-volume ones may not be. After removing two opposite facets from the 24-cell, we get an infinite-volume Coxeter polytope that can be deformed to form a continuous family of polytopes containing countably many Coxeter ones, few of them having finite-volume: this beautiful construction may be interpreted as a four-dimensional instance of the three-dimensional hyperbolic Dehn filling and is described in \cite{KSGT}.

\section{Hyperbolic four-manifolds} \label{manifolds:section}
We now list all the finite-volume hyperbolic four-manifolds that we were able to find in the literature. We introduce the subject with a paragraph that contains some general facts about finite-volume hyperbolic $n$-manifolds with $n\geq 4$. Many of them were proved in the last 15 years.

\subsection{General facts}
How can one construct hyperbolic manifolds in any dimension $n$? How many such manifolds are there, and how do they look like? Some of these questions have been answered for all $n$ simultaneously, mainly using arithmetic techniques. We now briefly state some theorems that hold in every fixed dimension $n\ge 4$, before turning back to dimension four in the following section.

We know from a theorem of Wang \cite{Wang} that there are only finitely many finite-volume complete hyperbolic $n$-manifolds of bounded volume. If we indicate by $\rho_n(V)$ the number of complete hyperbolic $n$-manifolds with volume $\leq V$, it was then shown by Burger, Gelander, Lubotzky, and Mozes in 2002 that $\rho_n(V)$ grows like $V^{cV}$, that is, there are two constants $0<c_1 < c_2$ such that
$$V^{c_1V} < \rho_n(V) < V^{c_2V}$$
for all sufficiently big $V$. The same kind of growth holds if one restricts to either closed or cusped manifolds, or to either arithmetic or non-arithmetic manifolds \cite{BGLS, BGLM}.

Gromov and Piatetski-Shapiro have first shown that there are non-arithmetic hyperbolic manifolds in all dimensions \cite{GPS}: these manifolds are constructed by gluing altogether two hyperbolic manifolds along their (isometric) geodesic boundaries. 

There are dramatically more non-arithmetic manifolds than arithmetic ones if one looks only at commensurability classes: Belolipetsky has shown \cite{Be2} that the number of commensurabillity classes of arithmetic hyperbolic $n$-manifolds with volume $\leq V$ grows roughly polynomially in $V$, while first Raimbault \cite{Rai}, and then Gelander and Levit \cite{GL}, have discovered that the number of commensurability classes of non-arithmetic ones grows first at least like $C^V$, and then like $V^{cV}$, respectively.

We also note that there are hyperbolic manifolds with arbitrarily small systole, as proved by Agol \cite{Ag} in dimension $n=4$ and then by Belolipetsky and Thomson \cite{BT} for every $n$.

Being reassured that there are plenty of hyperbolic four-manifolds, we now describe the (quite few, we must say) concrete examples that we have been able to find in the literature. We start by describing some closed hyperbolic four-manifolds, then we turn to finite-volume cusped ones, and finally to hyperbolic four-manifolds with (three-dimensional) geodesic boundary.

\subsection{The Davis manifold}
A couple of famous geometric three-manifolds are constructed by identifying the opposite faces of a regular dodecahedron $D$. Since opposite faces in $D$ are pentagons misaligned by a $\frac{\pi}5$ turn, to identify them we must make a choice: we can identifiy them through a (say, counterclockwise) rotation of angle $\frac{\pi}5$, $\frac{3\pi}5$, or $\pi$. If we choose the same angle on all opposite faces, we get a closed manifold $M$, which is correspondingly the Seifert-Weber manifold, Poincar\'e homology sphere, or $\matRP^3$. The orbit of every edge has order correspondingly 5, 3, 2: therefore if $D$ has dihedral angles $\frac{2\pi}5$, $\frac{2\pi}3$, $\pi$ the manifold $M$ is geometric. Such a dodecahedron $D$ exists in the appropriate geometry, and we get a hyperbolic manifold and two spherical manifolds respectively.

Davis made a similar (actually, simpler!) construction in 1985 using the 120-cell $P$, see \cite{Da}. Contrary to the dodecahedron, opposite facets in the 120-cell are parallel and not misaligned, hence one can just identify them via a translation. The orbit of a two-dimensional face has order 5 and this suggests that by realizing $P$ as the hyperbolic $\frac{2\pi}5$-angled polytope we get a hyperbolic manifold $M$. This is indeed the case, and $M$ is called the \emph{Davis manifold}.

Later on, Ratcliffe and Tschantz \cite{RT3} have studied $M$ more closely. They proved that $M$ is the unique smallest manifold among all coverings of the simplex orbifold $\Delta_3$. The degree of the covering is 14400, equal to the order of the symmetry group of the 120-cell. Its Euler characteristic is $\chi(M)=26$. The integral homology groups are as follows:
$$H_0(M) = \matZ, \quad H_1(M) = \matZ^{24}, \quad H_2(M) = \matZ^{72}, \quad H_3(M) = \matZ^{24}, \quad H_4(M) = \matZ.$$
The manifold $M$ contains as much as $720$ closed geodesics with the same minimal length $2.76514\ldots$ and therefore its injectivity radius is $1.38257\ldots$ This is the same injectivity radius of a closed hyperbolic surface $S$ obtained by identifying the opposite edges of a $\frac{2\pi}5$-angled hyperbolic regular decagon, and $M$ contains such a $S$. The isometry group of $M$ has order 28800.

Finally, the manifold $M$ is spin: since $H_1(M)$ has no torsion, this is in fact equivalent to the assertion that the intersection form of $H_2(M)$ is even. Evenness is proved by detecting 72 totally geodesic genus-2 surfaces that generate $H_2(M)$ and have trivial normal bundle: their self-intersection is zero, hence the form is even. 

Every closed oriented hyperbolic four-manifold has zero signature by the Hirzebruch signature formula (see \cite{LR1} for instance). Therefore the intersection form of a closed oriented hyperbolic four-manifold is either $\oplus_k \matr 0110$ or $\oplus_k (1) \oplus_k (-1)$ depending on whether the signature is even or odd. The intersection form of the Davis manifold is $\oplus_{36}\matr 0110$ and it seems to be the only closed four-manifold for which the signature and the intersection form have been determined: in particular, no closed oriented hyperbolic four-manifold with odd intersection form seems to be known.
Note that the second Betti number of an oriented closed hyperbolic four-manifold is even, because the Euler characteristic is so \cite{R}.

\subsection{The Conder--Maclachlan manifold}
For a certain time the Davis manifold $M$ was, with $\chi(M) = 26$, the smallest closed orientable hyperbolic four-manifold known. The situation changed in 2005 when Conder and Maclachlan made an extensive (and successful) search for a smaller closed orientable manifold. 

Their investigation went as follow: they picked the small compact Coxeter simplices $\Delta_1, \ldots, \Delta_5$ illustrated in Figure \ref{simplexes:fig}, and noted by some elementary considerations that the Euler characteristic of any orientable closed manifold $M$ covering these would be a multiple of (respectively) 2, 34, 26, 34, and 22. The Davis manifold has $\chi=26$ and is indeed a minimal orientable manifold covering of $\Delta_3$.

The only Coxeter simplex that can be covered by a closed orientable manifold smaller than the Davis manifold is hence $\Delta_1$. A computer search then allowed them to find a torsion-free subgroup $\Gamma < \Gamma_1$ of index $115200$. The quotient $N=\matH^4/_\Gamma$ is then a closed hyperbolic manifold with $\chi(N) = 115200\cdot \chi(\Delta_1) = \frac{115200}{14400} = 8$. The manifold $N$ is non-orientable, its orientable double cover $M$ has $\chi(M) = 16$ and is today the smallest closed orientable hyperbolic four-manifold known. Its homology groups are
$$
H_1(M) = \matZ^{2} \oplus \matZ/_{4\matZ} \oplus \left(\matZ/_{2\matZ}\right)^2, \quad
H_2(M) = \matZ^{18} \oplus \matZ/_{4\matZ} \oplus \left(\matZ/_{2\matZ}\right)^2, \quad
H_3(M) = \matZ^{2}.$$ 
The Betti numbers are smaller than those of the Davis manifold, but here homology has some torsion. The computer search was clever but not exhaustive, in the sense that $\Gamma_1$ could in principle contain some torsion-free subgroup of smaller index. 

In 2008 Long discovered more examples \cite{L} of torsion-free subgroups of $\Gamma_1$ of index $11520$ and hence of non-orientable closed hyperbolic four-manifolds with $\chi=8$.

\subsection{The Ratcliffe--Tschantz census}
The first systematic study of hyperbolic four-manifolds has been done by Ratcliffe and Tschantz in a 2000 paper \cite{RT0} that includes a census of $1171$ finite-volume hyperbolic cusped four-manifolds having minimal Euler characteristic $\chi = 1$.

In the hyperboloid model, the isometry group of $\matH^n$ is the group $O_+(n,1,\matR)$ of all positive lorentzian matrices, and we now look at the discrete subgroup $\Gamma = O_+(n,1,\matZ)$ consisting of all matrices with integer coefficients. In dimension $n=4$, Vinberg \cite{Vin67} proved in 1967 that $\Gamma$ is the reflection group of the Coxeter simplex $\Delta_7$ shown in Figure \ref{ideal_simplexes:fig}. As usual $\Gamma$ contains various finite-index torsion-free subgroups, and it is natural to consider the principal congruence-$k$ subgroup $\Gamma_k \triangleleft \Gamma$ that consists of all matrices that are congruent to the identity modulo $k$, for some $k\geq 2$.

The authors consider the case $k=2$. The group $\Gamma_2$ is not yet torsion-free: the authors show that $\Gamma_2$ is again a group generated by reflections, more precisely those of the right-angled Coxeter polytope $P_4$ defined in Section \ref{more:right:subsection}. 
Every torsion-free subgroup in $\Gamma_2$ has index at least $16$, and the authors show that there are precisely $1171$ of them having minimal index $16$. Each such subgroup has the ideal 24-cell $C$ as a fundamental domain (which is tessellated in $16$ copies of $P_4$) and hence gives rise to a hyperbolic four-manifold, obtained by pairing the faces of $C$.

Since $\chi(C)=1$ and $C$ is non-compact, the result of this investigation is a list of $1171$ cusped manifolds $M$ with minimum Euler characteristic $\chi(M)=1$ and hence minimum volume $\Vol(M)=\frac{4\pi^2}3$. Only 22 of these manifolds are orientable. 

\begin{figure}
\begin{center}
\includegraphics[width = 8 cm] {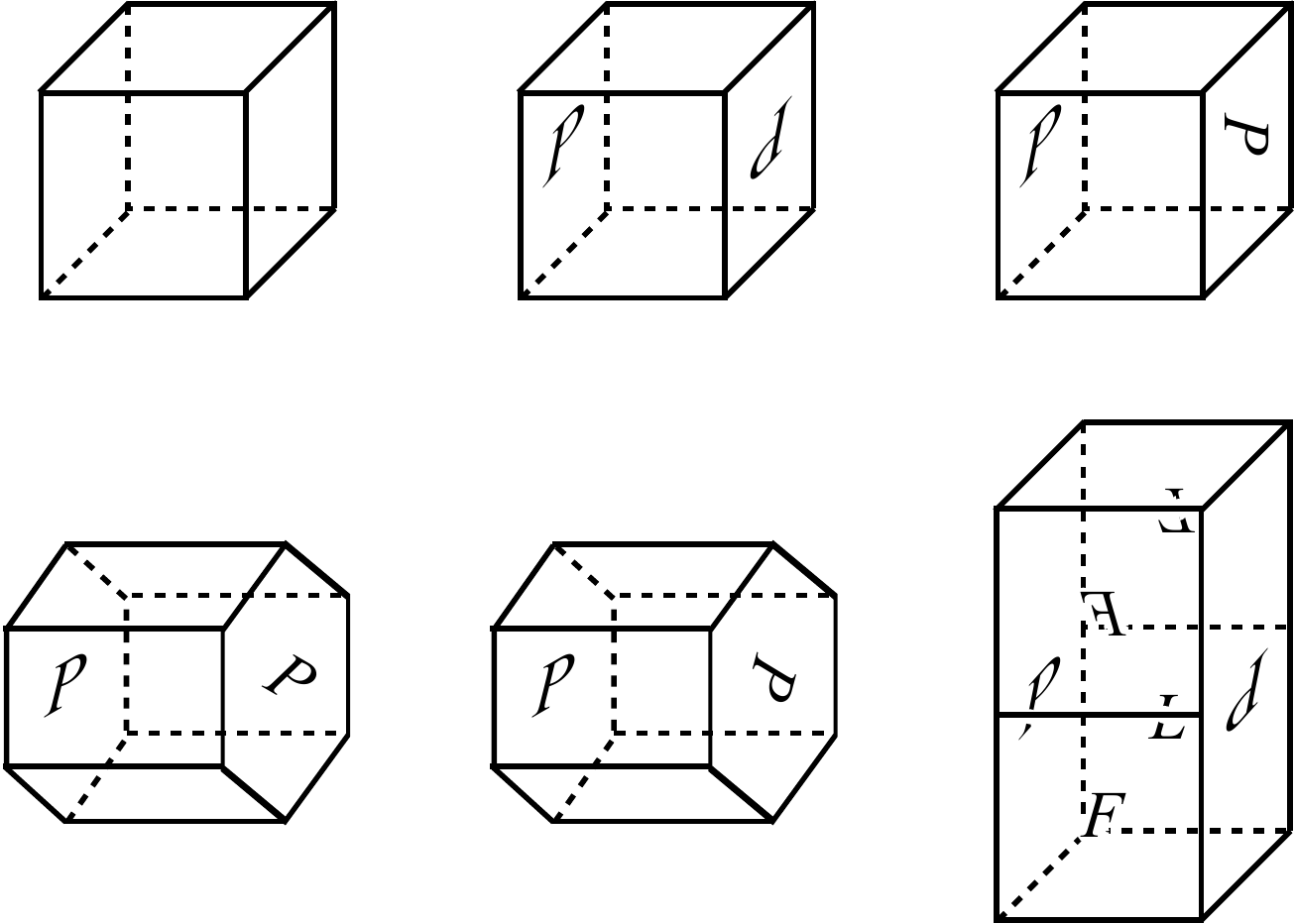}
\caption{The six closed orientable flat 3-manifolds, up to diffeomorphism. Each is constructed by pairing isometrically the faces of a polyhedron in $\matR^3$ according to the labels. When a face has no label, it is simply paired to its opposite by a translation. The polyhedra shown here are three cubes, two prisms with regular hexagonal basis, and one parallelepiped made of two cubes.}
\label{flat:fig}
\end{center}
\end{figure}

All the manifolds in the list have either 5 or 6 cusps: this is not surprising since they cover the orbifold $P_4$ that has $5$ ideal vertices. What are the cusp shapes of these manifolds? Recall that there are 6 orientable and 4 non-orientable closed flat 3-manifolds up to diffeomorphisms: the 6 orientable ones are shown in Figure \ref{flat:fig}, and they are Seifert manifolds over the orbifolds
$$T,\ (S^2,2,2,2,2),\ (S^2,2,4,4),\ (S^2,2,3,6),\ (S^2,3,3,3),\ (\matRP^2,2,2).$$
The first five manifolds in the list fiber over $S^1$ with torus fibers and monodromies
$$\begin{pmatrix} 1 & 0 \\ 0 & 1 \end{pmatrix}, \ 
\begin{pmatrix} -1 & 0 \\ 0 & -1 \end{pmatrix}, \ 
\begin{pmatrix} 0 & -1 \\ 1 & 0 \end{pmatrix}, \ 
\begin{pmatrix} 1 & -1 \\ 1 & 0 \end{pmatrix}, \ 
\begin{pmatrix} -1 & 1 \\ -1 & 0 \end{pmatrix}
$$
of order $1, 2, 4, 6$, and $3$ respectively. The sixth manifold is called the \emph{Hantzsche-Wendt manifold}, and it does not fiber over $S^1$ because its first homology group is finite.

Each of the $22$ orientable manifolds in the Ratcliffe--Tschantz census has 5 cusps: only the first two and the last among the six orientable diffeomorphism types appear as cusp sections of some of them. On the other hand, all the four non-orientable types occur among the $1171$ manifolds  \cite{RT0}.

It is worth saying that the $1171$ manifolds found are obtained from $C$ via a particular class of facet-pairings that uses only reflections along the coordinate hyperplanes. If one allows all kinds of isometries between facets, more manifolds are likely to be found: the number of hyperbolic manifolds obtained by pairing the facets of $C$ is completely unknown and can be much bigger than $1171$.

\subsection{Cusp shapes}
The Ratcliffe--Tschantz census shows that at least three among the six diffeomorphism types of orientable flat three-manifolds arise as cusp shapes of some finite-volume hyperbolic four-manifold. What about the other three?

Nimershiem proved in 1998 that each of the six orientable diffeomorphism types arises as a cusp shape of some multi-cusped finite-volume hyperbolic four-manifold \cite{N}. Much more than this, she proved that the geometric cusp shapes that arise in this way form a (countable) dense set in the (uncountable) flat moduli space of each of the six types.

Every diffeomorphism type appears as a cusp shape of a multi-cusped manifold, but there are some remarkable restrictions on the topology of all cusps considered altogether, discovered by Long and Reid \cite{LR1}. If $M$ is a cusped oriented hyperbolic four-manifold, its cusp shapes $N_1,\ldots, N_k$ are oriented flat three-manifolds, and the following equality holds:
\begin{equation} \label{eta:eqn}
\sigma (M) = \sum_{i=1}^k \eta(N_i)
\end{equation}
relating the signature $\sigma (M)$ of $M$ with the \emph{$\eta$-invariants} of the cusp shapes, see \cite{LR1}. The invariant $\eta(N)$ is a real number that depends only on the oriented diffeomorphism type of the flat manifold $N$, and changes sign on orientation reversals, hence in particular $\eta(N)=0$ when $N$ is mirrorable. The $\eta$-invariant of the six flat three-manifolds listed above were computed in \cite{Sz, O} and is (up to sign) equal to respectively
$$0, \ 0,\ 1,\ \frac 43,\ \frac 23,\ 0.$$

This shows in particular that the two flat manifolds with non-integral $\eta$-invariant cannot arise as the cusp shapes of a single-cusped hyperbolic four-manifold.

Concerning (\ref{eta:eqn}), it is worth noting that no cusped oriented hyperbolic manifold with non-zero signature seems to be known. (The signature of a \emph{closed} oriented hyperbolic manifold vanishes.)

\subsection{Manifolds with one cusp}
No single-cusped hyperbolic four-manifold at all was known until 2013, when many examples were built by Kolpakov and Martelli \cite{KM}. In their paper, they constructed infinitely many hyperbolic four-manifolds with any given number $k>0$ of cusps: the construction actually shows that the number $\rho_k(V)$ of $k$-cusped hyperbolic four-manifolds with volume bounded by $V$ grows like $V^{cV}$, for every $k>0$.

\begin{figure}
 \begin{center}
  \includegraphics[width = 3.5 cm]{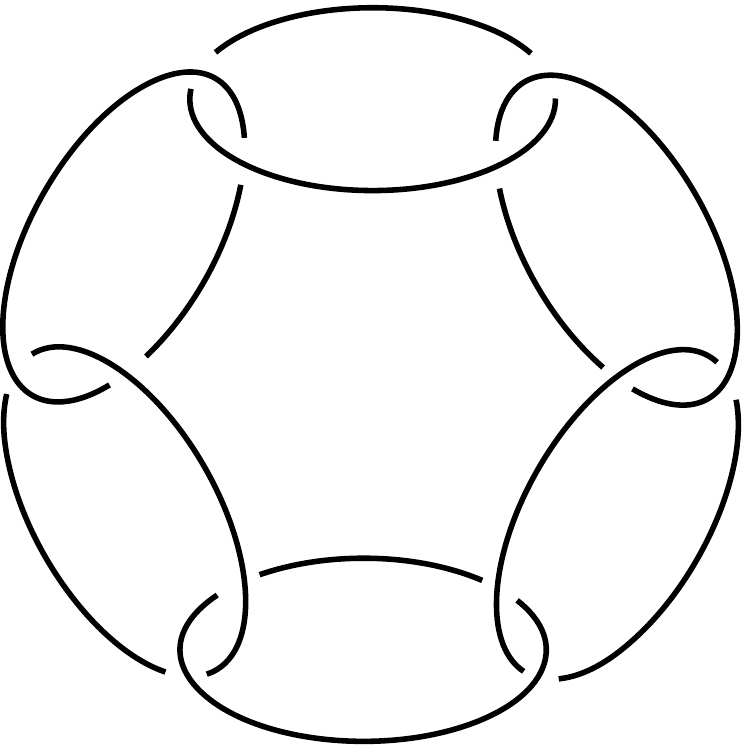}
 \end{center}
 \caption{The minimally twisted chain link with 6 components. Its complement $M$ is a hyperbolic three-manifold that is tessellated into four ideal regular octahedra. The block $B$ used in \cite{KM} has eight geodesic boundary components, each isometric to $M$.}  \label{chainlink:fig}
\end{figure}

All these manifolds are built by assembling multiple copies of a single block $B$, a finite-volume hyperbolic manifold with non-compact geodesic boundary and some rank-2 cusps. The block is constructed by noting that the 24 facets of the ideal 24-cell $C$ can be colored naturally in green, red, and blue, so that two adjacent facets have distinct colors. By doubling $C$ first along the green and then along the red strata we get $B$. The remaining blue strata now form the $8$ geodesic boundary components of $B$, each isometric to the complement of the chain link in $S^3$ with $6$ components shown in Figure \ref{chainlink:fig}. The block $B$ has $24$ rank-2 cusps of shape $S^1\times S^1\times [0,1]$, each connecting two distinct boundary components. One can assemble copies of $B$ in a way that all these rank-2 cusps glue up to form an arbitrary number $k$ of rank-3 cusps.

In all the single-cusped examples from \cite{KM}, the cusp shape is a 3-torus. More recently, in 2015 Kolpakov and Slavich \cite{KS2} constructed a single-cusped hyperbolic four-manifold with cusp shape the torus fibering with monodromy $\matr {-1}00{-1}$.

We still do not know if the torus fibering with monodromy $\matr 0{-1}10$ and the Hantzsche-Wendt manifold arise as cusp shapes of single-cusped hyperbolic four-manifolds.

\subsection{Geodesic boundary}
After constructing examples of closed and cusped hyperbolic four-manifolds, it is natural to consider hyperbolic manifolds with geodesic boundary. The first question that can arise is the following: are there finite-volume hyperbolic four-manifolds with connected geodesic boundary? The answer is affirmative, and the first compact and cusped examples were constructed by Ratcliffe and Tschantz in \cite{RT1} by suitably modifying the Davis and 24-cell constructions mentioned above. This question was also motivated by physical considerations \cite{G}.

Following Long and Reid \cite{LR1}, we say that a finite-volume orientable hyperbolic $(n-1)$-manifold $M$ \emph{bounds geometrically} if there is a finite-volume orientable hyperbolic $n$-manifold $W$ with geodesic boundary isometric to $M$. 

Which orientable hyperbolic three-manifolds bound geometrically? In 2000 Long and Reid proved \cite{LR1} that ``most'' hyperbolic closed three-manifolds do \emph{not} bound geometrically: indeed, the $\eta$-invariant $\eta(M)$ of a geometrically bounding closed hyperbolic three-manifold $M$ must be an integer, and this is a quite strong requirement because the surgeries on a hyperbolic knot provide a set of hyperbolic manifolds whose $\eta$-invariants form a dense subset of $\matR$, see \cite{MN}. 

In a 2001 paper \cite{LR3}, Long and Reid then proved that, despite this strong restriction, there are infinitely many commensurability classes of geometrically bounding hyperbolic manifolds in all dimensions. A natural quest is now to try to determine the smallest ones.

The smallest closed hyperbolic three-manifold $M$ with $\eta(M)\in\matZ$ is the arithmetic manifold named Vol3, with volume equal to the volume $1.0149 \ldots$ of the ideal regular tetrahedron \cite{LR1}. We still do not know whether this manifold bounds geometrically or not. It has $\eta(M) = 0$.

For the time being, the smallest closed geometrically bounding hyperbolic three-manifold known has volume $68.8992\ldots$ and was constructed by Kolpakov, Martelli, and Tschantz \cite{KMT} by assembling some 120-cells using a \emph{coloring} technique already employed by various authors in similar contexts, see for instance Vesnin \cite{Vesnin87, Vesnin}, Davis and Januszkiewicz \cite{DJ}, Izmestiev \cite{I}. 

\begin{figure}
 \begin{center}
  \includegraphics[width = 4 cm]{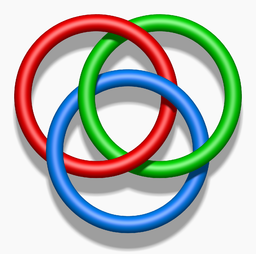}
  \includegraphics[width = 3.3 cm]{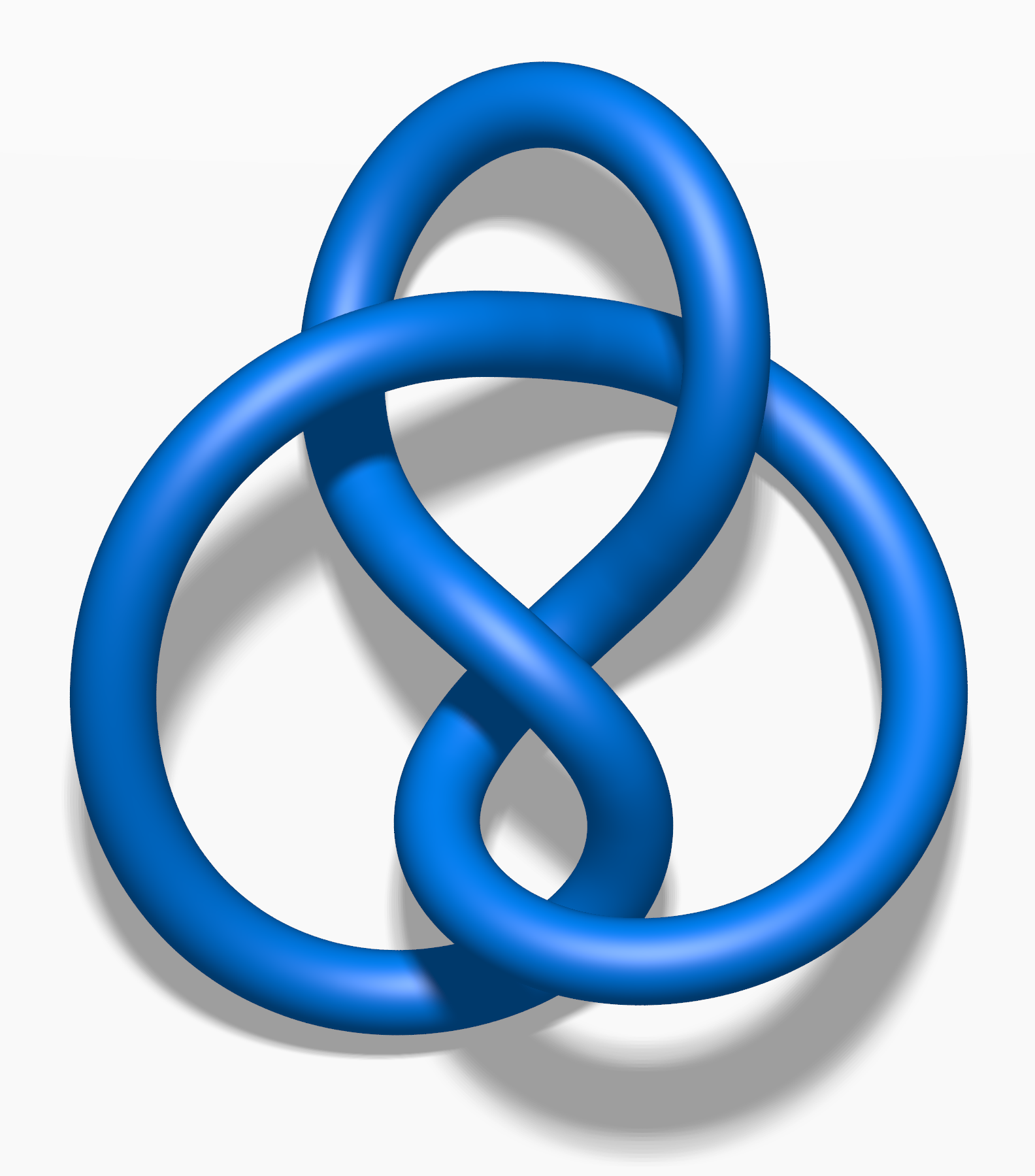}
 \end{center}
 \caption{The Borromean rings and figure-eigth knot complements are hyperbolic and are tessellated into two ideal regular octahedra and tetrahedra, respectively. They both bound geometrically some hyperbolic four-manifold that are tessellated into ideal 24-cells and rectified simplices, respectively.}  \label{Borromean:fig}
\end{figure}

On the cusped side, more progress have been done recently. The first cusped geometrically bounding 3-manifolds were constructed by Ratcliffe and Tschantz \cite{RT1}, then came an explicit link complement in the three-sphere \cite{S}, and very recently it was discovered that the Borromean link \cite{M} and the figure-eight knot complement \cite{S2} drawn in Figure \ref{Borromean:fig} both bound geometrically.

The figure-eight knot complement has volume 2.029 $\ldots$ and is the smallest orientable cusped hyperbolic three-manifold, together with its sibling \cite{CaMe}. Therefore the figure-eight knot complement is the smallest cusped geodesically bounding hyperbolic three-manifold. As shown by Slavich \cite{S2}, it bounds an orientable hyperbolic four-manifold $W$ with $\chi(W)=2$. The manifold $W$ is tessellated into 12 ideal rectified simplices.




\section{Manipulating hyperbolic four-manifolds} \label{manipulating:section}
In this section we describe some more elaborate constructions that involve hyperbolic four-manifolds. 

\subsection{Dehn filling}
Let $M$ be a finite-volume hyperbolic $n$-manifold with some $k\geq 1$ cusps, equipped with disjoint embedded cusp sections $X_1, \ldots, X_k$. Each section $X_i$ is a closed flat $(n-1)$-manifold, and we suppose that it is diffeomorphic to a $(n-1)$-torus: this not a too restrictive requirement, because it is always fulfilled virtually (ie on some finite cover), as a consequence of residual finiteness of $\pi_1(M)$ and the Bieberbach Theorem.

Every primitive element $\gamma_i \in H_1(X_i,\matZ)$ is represented by a simple closed geodesic, unique up to translations. The \emph{Dehn filling} of $M$ that \emph{kills} $\gamma_i$ is the topological operation that consists of truncating $M$ along $X_i$ and then shrinking all closed geodesic in $X_i$ parallel to $\gamma_i$ to points. Topologically, this is like gluing the ``solid torus'' $N=D^2 \times (S^1\times \ldots \times S^1)$ to the truncated $M$ via a diffeomorphism $\psi\colon \partial N \to X_i$ that sends $\partial D^2 \times \{{\rm pt}\}$ to $\gamma_i$.

A full Dehn filling of $M$ is parametrized by some primitive elements $(\gamma_1,\ldots, \gamma_k)$, one for each cusp. The result of this operation is a closed smooth four-manifold $M^{\rm fill}$.
In dimension $n=3$, the celebrated Thurston Dehn Filling Theorem states that, if one avoids finitely many ``exceptional'' classes $\gamma_i$ at each $X_i$, the filled manifold $M^{\rm fill}$ is again hyperbolic, with a hyperbolic metric that is obtained by deforming globally the initial one on $M$. This phenomenon is intrinsically three-dimensional and cannot be extended as it is to higher dimension $n\geq 4$. However, the theorem can be extended in a couple of interesting ways.

The first is the Gromov-Thurston $2\pi$ Theorem, which says that if every closed geodesic $\gamma_i$ has length $\ell(\gamma_i)\geq 2\pi$ then $M^{\rm fill}$ admits a metric of non-positive sectional curvature: the metric is obtained by extending the original hyperbolic one on the truncated $M$ to the attaching solid tori $N$, see \cite{A} for a proof.

The second is a more difficult theorem, proved by Anderson \cite{A} and Bamler \cite{B} via analytic methods: if each $\gamma_i$ avoids finitely many exceptional elements in $H_1(X_i, \matZ)$, the manifold $M^{\rm fill}$ admits an Einstein metric with negative scalar curvature. Note that such a metric is necessarily hyperbolic in dimension $n=3$ up to rescaling, but when $n\geq 4$ its sectional curvature is not necessarily non-positive everywhere (although it is negative in the average).

The Gromov-Thurston $2\pi$ Theorem, together with the Cartan-Hadamard Theorem, imply that if $\ell(\gamma_i)\geq 2\pi$ for all $i$ then $M^{\rm fill}$ is covered by $\matR^n$ and is in particular aspherical. 

In dimension four, we note that 
\begin{equation} \label{fill:eqn}
\chi(M^{\rm fill}) = \chi(M) >0, \quad \sigma(M^{\rm fill}) = \sigma(M) = 0, \quad \|M^{\rm fill}\| \leq \|M\|.
\end{equation}
Here $\sigma$ and $\|\cdot \|$ denote signature and simplicial volume (the first is defined only when $M$ is oriented). The two signatures are equal by the Novikov additivity Theorem, and $\sigma(M)=0$ because of (\ref{eta:eqn}): the $\eta$-invariant of the boundary is zero because it consists of 3-tori. The inequality on simplicial volumes is a theorem of Fujiwara and Manning \cite{FM} and it is still unknown whether it can be promoted to a strict inequality, as in dimension three.

\subsection{Aspherical homology spheres}
The Gromov-Thurston $2\pi$ Theorem was used by Ratcliffe and Tschantz to construct the first examples of closed aspherical integral homology spheres in dimension four \cite{RT4}. 

Their construction goes as follows: they pick a particularly symmetric non-orientable manifold $M$ in their census, that has $\chi(M)=1$ and five cusps, and consider its orientable double covering $\tilde M$, which has $\chi(\tilde M) = 2$ and five 3-torus cusps. They compute that $H_1(\tilde M) = \matZ^5$ and show that for many choices of $(\gamma_1,\ldots, \gamma_5)$ the Dehn filling kills all the first homology and hence produces a $\tilde M^{\rm fill}$ with $\chi = 2$ and $H_1=H_3=\{e\}$. Therefore also $H_2=\{e\}$ and $\tilde M^{\rm fill}$ is a homology sphere. One can choose $(\gamma_1,\ldots,\gamma_5)$ with arbitrarily big $\ell(\gamma_i)$, hence if $\ell(\gamma_i)>2\pi $ the manifold $\tilde M^{\rm fill}$ is aspherical.

It seems still unknown whether there is a closed homology sphere $M$ that admits a negative sectional curvature metric, or even a hyperbolic one: such a hyperbolic $M$ would have $\chi(M)=2$ and hence have minimum volume, since closed orientable hyperbolic manifolds have even Euler characteristic \cite{R}.

\subsection{Link complements}
Let $M$ be a cusped hyperbolic four-manifold as above, with cusp sections $X_1,\ldots, X_k$ all diffeomorphic to 3-tori. If $M^{\rm fill}$ is a Dehn filling of $M$, the cores of the attached $k$ solid tori form $k$ disjoint smooth two-dimensional tori in $M^{\rm fill}$. One can hence represent $M$ as the complement of a link of $k$ tori in $M^{\rm fill}$, similarly as one represents a cusped hyperbolic 3-manifold as a link complement in any of its Dehn fillings. 

As in three dimensions, let us say that a set of $k$ linked tori in some closed four-manifold $W$ is \emph{hyperbolic} if its complement admits a complete hyperbolic metric. In dimension three, every closed three-manifold contains plenty of hyperbolic links. What happens in dimension four?

It is natural to expect that most closed four-manifolds $W$ contain no hyperbolic link at all. This is certainly the case if $\chi(W)\leq 0$ or $\sigma(W)\neq 0$, see (\ref{fill:eqn}). Moreover, by Wang's Theorem \cite{Wang} there are only finitely many cusped $M$ with any given Euler characteristic. There are infinitely many Dehn fillings of such manifolds $M$, but they are likely to cover only a tiny portion of the set of all closed four-manifolds with that Euler characteristic. They certainly have bounded first homology groups, for instance.

Despite these pessimistic premises, some hyperbolic link complements in familiar closed four-manifolds have been found in the last years. Ivan\v{s}i\'c \cite{Iv} and then Ivan\v{s}i\'c, Ratcliffe, and Tschantz \cite{IRT} have constructed twelve distinct manifolds that are link complement in some homotopy 4-sphere (the links contain tori and/or Klein bottles). Their construction goes as follows: they look at the orientable double covers of the $1171$ Ratcliffe-Tschantz manifolds, that have $\chi =2$, and check whether some of them could have a simply connected Dehn filling. If this is the case, the filled manifold is a homotopy sphere because $\chi=2$. The authors could find such a Dehn filling in 12 cases. These homotopy 4-spheres are homeomorphic to $S^4$ thanks to Freedman's Theorem, but might a priori be not diffeomorphic to it.

Later on, Ivan\v{s}i\'c proved \cite{Iv2} that one of these homotopy spheres is indeed diffeomorphic to $S^4$. We therefore know that $S^4$ contains a hyperbolic link of 5 tori. More recently, Saratchandran proved \cite{Sa, Sa2} that there are hyperbolic link complements in $S^2\times S^2$ and in some four-manifolds homeomorphic to $\#_{2k}(S^2\times S^2)$ for every $k\geq 0$. All these results are proved by Dehn-filling some finite covers of manifolds from the Ratcliffe--Tschantz census.

\subsection{Symmetries}
A finite-volume hyperbolic manifold has finite isometry group. Conversely, Belolipetsky and Lubotzky have shown \cite{BL} that every finite group $G$ arises as the isometry group of a $n$-dimensional finite-volume hyperbolic manifold $M$, for all $n\geq 2$.

In 2015 Kolpakov and Slavich \cite{KS} reproved this theorem in dimension four, producing for every $G$ a hyperbolic four-manifold $M$ with $\Vol(M) \leq C |G|\log^2 |G|$ for some fixed constant $C>0$. They also showed that the number of manifolds with fixed isometry group $G$ and volume smaller than $V$ grows like $V^{cV}$.

The manifold $M$ is constructed by assembling multiple copies of the same block with geodesic boundary. The block is built using the ideal rectified simplex described in Section \ref{rectified:subsection}.

\subsection{Symplectic structures}
How do gauge invariants behave on hyperbolic four-manifolds? Do some closed hyperbolic four-manifolds admit some nice structure, like a symplectic structure or a Lefschetz fibration? 

These questions are still wide open, however there is a pessimistic conjecture around in the literature: it is conjectured by LeBrun \cite{LB} that the Seiberg-Witten invariants of every closed hyperbolic four-manifold should vanish. Note that, by a celebrated theorem of Taubes \cite{T}, a closed symplectic four-manifold with $b_2^+ \geq 2$ has non-zero Seiberg-Witten invariants, and in addition every closed four-manifold admitting a Lefschetz fibration with a homologically non-trivial fiber has a symplectic structure \cite[Theorem 10.2.18]{GS}. So both structures (symplectic and Lefschetz fibration) seem unlikely to exist on closed hyperbolic four-manifolds, at least following the LeBrun conjecture.

On the other hand, closed hyperbolic four-manifolds like the Davis manifold can be used to construct symplectic 6-manifolds \cite{FP}.

\section{Some open questions} \label{questions:section}
It is relatively easy to formulate open questions on hyperbolic four-manifolds, since only few concrete examples are known. We have already alluded to some open problems in the previous pages: we collect some of them here, and we add a few more. As in the rest of the paper, all hyperbolic manifolds are assumed to be complete and finite-volume.

\begin{enumerate}
\item[1.] Can we find a closed oriented hyperbolic four-manifold with odd intersection form? 
\item [2.] Can we find a cusped oriented hyperbolic four-manifold with non-vanishing signature? What are the intersection forms of cusped oriented hyperbolic four-manifolds?
\item [3.] Are there infinitely many cusped or closed hyperbolic four-manifolds with bounded first Betti number? Or, more strongly, with a bounded number of generators for their fundamental groups?
\item[4.] Which groups arise as homology groups of hyperbolic four-manifolds?
\item [5.] What is the smallest volume of a closed oriented hyperbolic four-manifold?
\item [6.] How many cusped four-manifolds are obtained by coupling the facets of an ideal 24-cell? Is there a cusped four-manifold with $\chi=1$ that does not arise from this construction? Are all hyperbolic four-manifolds with $\chi=1$ commensurable?
\item [7.] Is there a hyperbolic four-manifold that fibers nicely in any way?
\item [8.] Is there a way to write Thurston's equations for ideal triangulations in dimension four? 
\item [9.] Does every hyperbolic four-manifold contain immersed geodesic hyperbolic surfaces and 3-manifolds?
\item [10.] Can we model numerically the Einstein metrics on the Dehn fillings of some cusped hyperbolic four-manifold?
\end{enumerate}

Concerning question 6, note that a hyperbolic four-manifold cannot fiber over $S^1$ because it has positive Euler characteristic. One could however look at some nicely behaved $ S^1$-valued Morse functions, or fibrations over surfaces with some controlled singularities. 


\end{document}